\documentclass[12pt]{article}
\usepackage{amsmath}

\begin{document}

\begin{titlepage}
\title{\bf Lagrangian and Hamiltonian Dynamics \\on Para-K\"{a}hlerian Space Form}
\author{Mehmet Tekkoyun \footnote{tekkoyun@pamukkale.edu.tr} \\
 {\small Department of Mathematics, Pamukkale University,}\\
{\small 20070 Denizli, Turkey}}
\date{\today}
\maketitle

\begin{abstract}

In this study, we
introduce Euler-Lagrange and Hamiltonian equations on \ $(%
\mathbf{R}%
^{2},g,J)$ being a model of para-K\"{a}hlerian Space Forms.
Finally,
some geometrical and physical results on the related mechanic
systems have
been discussed.



\end{abstract}
\end{titlepage}

\section{Introduction}

Modern Differential Geometry is a suitable frame for studying Lagrangian and
Hamiltonian formalisms of Classical Mechanics. To show this, it is possible
to find many articles and books in the relevant fields. It is well-known
that the dynamics of Lagrangian and Hamiltonian systems is characterized by
a convenient vector field $X$ defined on the tangent and cotangent bundles
which are phase-spaces of velocities and momentum of a given configuration
manifold. If $Q$ is an $m$-dimensional configuration manifold and $%
L:TQ\rightarrow R$ is a regular Lagrangian function, then there is a unique
vector field $X$ on $TQ$ such that
\begin{equation}
i_{X_{L}}\omega _{L}=dE_{L},  \label{1.1}
\end{equation}%
where $\omega _{L}$ is the symplectic form and $E_{L}$ is energy associated
to $L$. The so-called Euler-Lagrange vector field $X$ is a semispray (or
\textit{second order differential equation}) on $Q$ since its integral
curves are the solutions of the Euler-Lagrange equations. The triple, either
$(TQ,\omega _{L},\xi )$ or $(TQ,\omega _{L},L),$ is called \textit{%
Lagrangian system} on the tangent bundle $TQ.$ If $H:T^{\ast }Q\rightarrow R$
is a regular Hamiltonian function then there is a unique vector field $X_{H}$
on $T^{\ast }Q$ such that
\begin{equation}
i_{X_{H}}\omega =dH  \label{1.2}
\end{equation}%
where $\omega $ is the symplectic form and $H$ stands for Hamiltonian
function. The paths of the so-called Hamiltonian vector field $X_{H}$ are
the solutions of the Hamiltonian equations. The triple, either $(T^{\ast
}Q,\omega ,Z_{H})$ or $(T^{\ast }Q,\omega ,H),$ is called \textit{%
Hamiltonian system} on the cotangent bundle $T^{\ast }Q$ fixed with
symplectic form $\omega $.

From the before some studies given in \cite{crampin, deleon, crampin1,
tekkoyun1, tekkoyun2, tekkoyun3, tekkoyun4}; we know that time-dependent or
not real, complex and paracomplex analogues of the Euler- Lagrange and
Hamiltonian equations have detailed been introduced. But, we see that is not
mentioned about Lagrangian and Hamiltonian dynamics on para-K\"{a}hlerian
space forms. Therefore, in this paper we present the Euler-Lagrange
equations and Hamiltonian equations on a model of para-K\"{a}hlerian space
forms and to derive geometrical and physical results on related dynamics
systems.

In this study, all the manifolds and geometric objects are $C^{\infty }$ and
the Einstein summation convention is in use. Also, $\mathbf{R}$, $\mathcal{F}%
(M)$, $\chi (M)$ and $\Lambda ^{1}(M)$ denote the set of real numbers, the
set of functions on $M$, the set of vector fields on $M$ and the set of
1-forms on $M$, respectively.

\section{Para-K\"{a}hlerian Space Forms}

\textbf{Definition 1 }\cite{bejan, cruceanu}\textbf{\ : }Let a manifold $M$
be endowed with an almost product structure $J\neq \mp Id$; which is a (1;
1)-tensor field such that $J^{2}=Id$: We say that $(M,J)$ (resp.$(M,J,g)$)
is an almost product (resp. almost Hermitian) manifold, where $g$ is a
semi-Riemannian metric on $M$ with respect to which $J$ is skew-symmetric,
that is
\begin{equation}
g(JX,Y)+g(X,JY)=0,\forall X,Y\in \chi (M)  \label{2.1}
\end{equation}%
Then $(M,J,g)$ is para-K\"{a}hlerian if $J$ is parallel with respect to the
Levi-Civita connection.

Let $(M,J,g)$ be a para-K\"{a}hlerian manifold and let denote the curvature
(0, 4)-tensor field by
\begin{equation}
R(X,Y,Z,V)=g(R(X,Y)Z,V);\forall X,Y,Z,V\in \chi (M)  \label{2.2}
\end{equation}%
where the Riemannian curvature (1, 3)-tensor field associated to the
Levi-Civita connection $\nabla $ of $g$ is given by $R=$ [$\nabla ,\nabla $]
-$\nabla _{\left[ \text{ },~\text{\ }\right] }$.

Then

\begin{equation}
\begin{array}{c}
R(X,Y,Z,V)=-R(Y,X,Z,V)=-R(X,Y,V,Z)=R(JX,JY,Z,V) \\
\text{and}\underset{\sigma }{\text{ }\sum }R(X,Y,Z,V)=0,%
\end{array}
\label{2.3}
\end{equation}

where $\sigma $ denotes the sum over all cyclic permutations. We know that
the following (0,4)-tensor field is defined by%
\begin{equation}
R_{0}(X,Y,Z,V)=\frac{1}{4}\left\{
\begin{array}{c}
g(X,Z)g(Y,V)-g(X,V)g(Y,Z)-g(X,JZ)g(Y,JV) \\
+g(X,JV)g(Y,JZ)-2g(X,JY)g(Z,JV)%
\end{array}%
\right\} ,  \label{2.4}
\end{equation}

where $\forall X,Y,Z,V\in \chi (M).$ For any $p\in M$, a subspace $S\subset
T_{p}M$ is called non-degenerate if $g$ restricted to $S$ is non-degenerate.
If $\left\{ u,v\right\} $ is a basis of a plane $\sigma \subset T_{p}M$,
then $\sigma $ is

non-degenerate iff $g(u,u)g(v,v)-[g(u,v)]^{2}\neq 0$. In this case the
sectional curvature of $\sigma $= span$\left\{ u,v\right\} $ is

\begin{equation}
k(\sigma )=\frac{R(u,v,u,v)}{g(u,u)g(v,v)-[g(u,v)]^{2}}  \label{2.5}
\end{equation}

From (\ref{2.1}) it follows that $X$ and $JX$ are orthogonal for any $X\in $
$\chi (M)$. By a $J$-plane we mean a plane which is invariant by $J$. For
any $p\in M$, a vector $u$ $\in T_{p}M$ is isotropic provided $g(u,u)=0$. If
$u$ $\in T_{p}M$ is not isotropic, then the sectional curvature $H(u)$ of
the $J$-plane span$\left\{ u,Ju\right\} $ is called the $J$-sectional
curvature defined by $u.$ When $H(u)$ is constant,

then $(M,J,g)$ is called of constant $J$-sectional curvature, or a para-K%
\"{a}hlerian space form.

\textbf{Theorem 1:} Let $(M,J,g)$ be a para-K\"{a}hlerian manifold such that
for each $p\in M$, there exists $c_{p}\in R$ satisfying $H(u)=c_{p}$ for $u$
$\in T_{p}M$ \ such that $g(u,u)g(Ju,Ju)\neq 0.$Then the Riemann-
Christoffel tensor $R$ satisfies $R=cR_{0},$ where $c$ is the function
defined by $p\rightarrow c_{p}.$ And conversely.

\textbf{Definition 2: }A para-K\"{a}hlerian manifold $(M,J,g)$ is said to be
of constant paraholomorphic sectional curvature $c$ if it satisfies the
conditions of \textbf{Theorem 1}.

\textbf{Theorem 2:} Let $(M,J,g)$ be a para-K\"{a}hlerian manifold with $%
dimM>2$. Then the following properties are equivalent:

1) $M$ is a space of constant paraholomorphic sectional curvature $c$

2) The Riemann- Christoffel tensor curvature tensor $R$ has the expression%
\begin{equation}
R(X,Y,Z,V)=\frac{c}{4}\left\{
\begin{array}{c}
g(X,Z)g(Y,V)-g(X,V)g(Y,Z)-g(X,JZ)g(Y,JV) \\
+g(X,JV)g(Y,JZ)-2g(X,JY)g(Z,JV)%
\end{array}%
\right\} ,  \label{2.6}
\end{equation}

where $\forall X,Y,Z,V\in \chi (M).$ Let $(x,\ y)$ be a real coordinate
system on a neighborhood $U$ of any point $p$ of $\mathbf{R}^{2},$ and $\{(%
\frac{\partial }{\partial x})_{p},(\frac{\partial }{\partial y})_{p}\}$ and $%
\{(dx)_{p},(dy)_{p}\}$ natural bases over $\mathbf{R}$ of the tangent space $%
T_{p}(\mathbf{R}^{2})$ and the cotangent space $T_{p}^{\ast }(\mathbf{R}%
^{2}) $ of $\mathbf{R}^{2},$ respectively.

The space $(\mathbf{R}^{2},g,J),$ is the model of the para-K\"{a}hlerian
space forms of dimension 2 and paraholomorphic sectional curvature $c\neq 0,$
where $g$ is the metric

\begin{equation}
g=\frac{4}{c}\left( \cosh ^{2}2ydx\otimes dx-dy\otimes dy\right) ,0\neq c\in
\mathbf{R,}  \label{2.7}
\end{equation}

and $J$ the almost product structure%
\begin{equation}
J=-\frac{1}{\cosh 2y}\frac{\partial }{\partial x}\otimes dy-\cosh 2y\frac{%
\partial }{\partial y}\otimes dx.  \label{2.8}
\end{equation}

Then we have
\begin{equation}
J(\frac{\partial }{\partial x})=-\cosh 2y\frac{\partial }{\partial y},\text{
}J(\frac{\partial }{\partial y})=-\frac{1}{\cosh 2y}\frac{\partial }{%
\partial x}.  \label{2.9}
\end{equation}%
The dual endomorphism $J^{\ast }$ of the cotangent space $T_{p}^{\ast }(%
\mathbf{R}^{2})$ at any point $p$ of manifold $\mathbf{R}^{2}$ satisfies $%
J^{\ast 2}=Id$ and is defined by
\begin{equation}
J^{\ast }(dx)=-\cosh 2ydy,\text{ }J^{\ast }(dy)=-\frac{1}{\cosh 2y}dx.
\label{2.10}
\end{equation}

\section{Lagrangian Dynamics}

Here, we find Euler-Lagrange equations for Classical Mechanics constructed
on para-K\"{a}hlerian space form $(\mathbf{R}^{2},g,J)$.

Denote by $J$ almost product structure and by $(x,y)$ the coordinates of $%
\mathbf{R}^{2}$. Assume that semispray be a vector field as follows:
\begin{equation}
\xi =X\frac{\partial }{\partial x}+Y\frac{\partial }{\partial y},\text{ }X=%
\overset{.}{x}=y,\text{ }Y=\overset{.}{y}.  \label{3.1}
\end{equation}
By \textit{Liouville vector field} on para-K\"{a}hlerian space form $(%
\mathbf{R}^{2},g,J),$ we call the vector field determined by $V=J\xi $ and
calculated by
\begin{equation}
J\xi =-\frac{1}{\cosh 2y}.Y\frac{\partial }{\partial x}-\cosh 2y.X\frac{%
\partial }{\partial y},  \label{3.2}
\end{equation}%
Given $T$ by \textit{the kinetic energy} and $P$ by \textit{the potential
energy of mechanics system} on para-K\"{a}hlerian space form. Then we write
by $L=T-P$ \textit{Lagrangian function }and by $E_{L}=V(L)-L$ \textit{the
energy function} associated $L$.

Operator $i_{J}$\ defined by\textbf{\ }%
\begin{equation}
i_{J}:\wedge ^{2}\mathbf{R}^{2}\rightarrow \wedge ^{1}\mathbf{R}^{2}\text{, }%
i_{J}(\omega )(X)=\omega (X,JX)  \label{3.3}
\end{equation}%
is called the \textit{interior product} with $J$, or sometimes the \textit{%
insertion operator}, or \textit{contraction} by $J$, where $\omega \in
\wedge ^{2}\mathbf{R}^{2},$ $X$ $\in \chi (\mathbf{R}^{2}).$ The exterior
vertical derivation $d_{J}$ is defined by
\begin{equation}
d_{J}=[i_{J},d]=i_{J}d-di_{J},  \label{3.4}
\end{equation}%
where $d$ is the usual exterior derivation. For almost product structure $J$
determined by (\ref{2.9}), the closed para-K\"{a}hlerian form is the closed
2-form given by $\Phi _{L}=-dd_{J}L$ such that
\begin{equation}
d_{J}=-\cosh 2y.\frac{\partial }{\partial y}dx-\frac{1}{\cosh 2y}.\frac{%
\partial }{\partial x}dy:\mathcal{F}(\mathbf{R}^{2})\rightarrow \wedge ^{1}%
\mathbf{R}^{2}.  \label{3.5}
\end{equation}%
Thus we get
\begin{eqnarray}
\Phi _{L} &=&\cosh 2y\frac{\partial ^{2}L}{\partial a\partial y}da\wedge
dx+\cosh 2y\frac{\partial ^{2}L}{\partial b\partial y}db\wedge dx  \notag \\
&&+\frac{1}{\cosh 2y}\frac{\partial ^{2}L}{\partial a\partial x}da\wedge dy+%
\frac{1}{\cosh 2y}\frac{\partial ^{2}L}{\partial b\partial x}db\wedge dy.
\label{3.6}
\end{eqnarray}%
where $(a,b)$ is other coordinates of $\mathbf{R}^{2}.$ Then%
\begin{equation}
\begin{array}{c}
i_{\xi }\Phi _{L}=\cosh 2y.X\frac{\partial ^{2}L}{\partial a\partial y}%
\delta _{a}^{x}dx-\cosh 2y.X\frac{\partial ^{2}L}{\partial a\partial y}%
da+\cosh 2y.Y\frac{\partial ^{2}L}{\partial b\partial y}\delta _{b}^{y}dx \\
-\cosh 2y.X\frac{\partial ^{2}L}{\partial b\partial y}db+\frac{1}{\cosh 2y}.X%
\frac{\partial ^{2}L}{\partial a\partial x}\delta _{a}^{x}dy-\frac{1}{\cosh
2y}.Y\frac{\partial ^{2}L}{\partial a\partial x}da \\
+\frac{1}{\cosh 2y}.Y\frac{\partial ^{2}L}{\partial b\partial x}\delta
_{b}^{y}dy-\frac{1}{\cosh 2y}.Y\frac{\partial ^{2}L}{\partial b\partial x}db.%
\end{array}
\label{3.7}
\end{equation}%
Since the closed para-K\"{a}hlerian form $\Phi _{L}$ on para-K\"{a}hlerian
space form $(\mathbf{R}^{2},g,J)$ is para-symplectic structure, one may find
\begin{equation}
E_{L}=\mathbf{-}\frac{1}{\cosh 2y}.Y\frac{\partial L}{\partial x}+\cosh 2y.X%
\frac{\partial L}{\partial y}-L,  \label{3.8}
\end{equation}%
and thus
\begin{equation}
\begin{array}{ll}
dE_{L}= & -\frac{1}{\cosh 2y}.Y\frac{\partial ^{2}L}{\partial a\partial x}%
da-\cosh 2y.X\frac{\partial ^{2}L}{\partial a\partial y}da-\frac{\partial L}{%
\partial a}da \\
& -\frac{1}{\cosh 2y}.Y\frac{\partial ^{2}L}{\partial b\partial x}db-\cosh
2y.X\frac{\partial ^{2}L}{\partial b\partial y}db-\frac{\partial L}{\partial
b}db.%
\end{array}
\label{3.9}
\end{equation}%
Considering $i_{\xi }\Phi _{L}=dE_{L}$, we calculate
\begin{equation}
\begin{array}{l}
\cosh 2y.X\frac{\partial ^{2}L}{\partial a\partial y}dx+\cosh 2y.Y\frac{%
\partial ^{2}L}{\partial b\partial y}dx \\
+\frac{1}{\cosh 2y}.X\frac{\partial ^{2}L}{\partial a\partial x}dy+\frac{1}{%
\cosh 2y}.Y\frac{\partial ^{2}L}{\partial b\partial x}dy+\frac{\partial L}{%
\partial x}dx+\frac{\partial L}{\partial y}dy=0.%
\end{array}
\label{3.10}
\end{equation}%
If the curve $\alpha :\mathbf{I\subset R}\rightarrow \mathbf{R}^{2}$ be
integral curve of $\xi ,$ which satisfies
\begin{equation}
\begin{array}{l}
\cosh 2y\left[ X\frac{\partial ^{2}L}{\partial a\partial y}+Y\frac{\partial
^{2}L}{\partial b\partial y}\right] dx+\frac{\partial L}{\partial x}dx \\
+\frac{1}{\cosh 2y}\left[ X\frac{\partial ^{2}L}{\partial a\partial x}+Y%
\frac{\partial ^{2}L}{\partial b\partial x}\right] dy+\frac{\partial L}{%
\partial y}dy=0,%
\end{array}
\label{3.11}
\end{equation}%
we get equations

\begin{equation}
\cosh 2y\frac{\partial }{\partial t}\left( \frac{\partial L}{\partial y}%
\right) +\frac{\partial L}{\partial x}=0,\frac{1}{\cosh 2y}\ \frac{\partial
}{\partial t}\left( \frac{\partial L}{\partial x}\right) +\frac{\partial L}{%
\partial y}=0  \label{3.12}
\end{equation}%
so-called \textit{Euler-Lagrange equations }whose solutions are the paths of
the semispray $\xi $ on para-K\"{a}hlerian space form $(\mathbf{R}^{2},g,J)$%
. Finally one may say that the triple $(\mathbf{R}^{2},\Phi _{L},\xi )$ is%
\textit{\ mechanical system} on para-K\"{a}hlerian space form $(\mathbf{R}%
^{2},g,J).$Therefore we say the following:

\textbf{Proposition 1: }Let $J$ almost product structure on para-K\"{a}%
hlerian space form $(\mathbf{R}^{2},g,J).$Also let $(f_{1},f_{2})$ be
natural bases of $\mathbf{R}^{2}.$Then it follows

\begin{equation*}
\begin{array}{cc}
\cosh 2y.J(f_{2})+f_{1}=0\Longleftrightarrow & \cosh 2y.\overset{.}{f}%
_{2,L}+f_{1,L}=0, \\
\frac{1}{\cosh 2y}J(f_{1})+f_{2}=0\Longleftrightarrow & \frac{1}{\cosh 2y}.%
\overset{.}{f}_{1,L}+f_{2,L}=0,%
\end{array}%
\end{equation*}

where $f_{1,L}=\frac{\partial L}{\partial x},$ $\ f_{2,L}=\frac{\partial L}{%
\partial y},$ $\overset{.}{f}_{1,L}=\frac{\partial }{\partial t}(\frac{%
\partial L}{\partial x}),$ $\ \overset{.}{f}_{2,L}=\frac{\partial }{\partial
t}(\frac{\partial L}{\partial y}).$

\section{Hamiltonian Dynamics}

Now, we conclude Hamiltonian equations for Classical Mechanics structured on
para-K\"{a}hlerian space form $(\mathbf{R}^{2},g,J)$.

Let $J^{\ast }$ be an almost product structure defined by (\ref{2.10}) and $%
\lambda $ Liouville form determined by $J^{\ast }(\omega )=-x\cosh 2ydy-y%
\frac{1}{\cosh 2y}dx$ such that $\omega =xdx+ydy$ 1-form on $\mathbf{R}^{2}.$
If $\Phi =-d\lambda $ is closed para-K\"{a}hlerian form$,$ then it is also a
para-symplectic structure on $\mathbf{R}^{2}$.

Let $(\mathbf{R}^{2},g,J)$\textbf{\ }be para-K\"{a}hlerian space form fixed
with closed para-K\"{a}hlerian form $\Phi $. Suppose that Hamiltonian vector
field $Z_{H}$ associated to Hamiltonian energy $H$ is given by
\begin{equation}
Z_{H}=X\frac{\partial }{\partial x}+Y\frac{\partial }{\partial y}.
\label{4.1}
\end{equation}

For the closed para-K\"{a}hlerian form $\Phi $ on $\mathbf{R}^{2},$ we have
\begin{equation}
\Phi =-d\lambda =-d(-x\cosh 2ydy-y\frac{1}{\cosh 2y}dx)=\frac{\cosh ^{2}2y-1%
}{\cosh 2y}dx\wedge dy.  \label{4.2}
\end{equation}%
Then it follows
\begin{equation}
i_{Z_{H}}\Phi =i_{Z_{H}}(-d\lambda )=-\frac{\cosh ^{2}2y-1}{\cosh 2y}Ydx+%
\frac{\cosh ^{2}2y-1}{\cosh 2y}Xdy.  \label{4.3}
\end{equation}%
Otherwise, we find the differential of Hamiltonian energy the following as
\begin{equation}
dH=\frac{\partial H}{\partial x}dx+\frac{\partial H}{\partial y}dy.
\label{4.4}
\end{equation}%
From (\ref{4.3}) and (\ref{4.4}) with respect to $i_{Z_{H}}\Phi =dH,$ we
find para-Hamiltonian vector field on para-K\"{a}hlerian space form to be
\begin{equation}
Z_{H}=\frac{\cosh 2y}{\cosh ^{2}2y-1}\frac{\partial H}{\partial y}\frac{%
\partial }{\partial x}-\frac{\cosh 2y}{\cosh ^{2}2y-1}\frac{\partial H}{%
\partial x}\frac{\partial }{\partial y}.  \label{4.5}
\end{equation}

Assume that the curve
\begin{equation}
\alpha :I\subset \mathbf{R}\rightarrow \mathbf{R}^{2}  \label{4.6}
\end{equation}%
be an integral curve of Hamiltonian vector field $Z_{H},$ i.e.,
\begin{equation}
Z_{H}(\alpha (t))=\overset{.}{\alpha },\,\,t\in I.  \label{4.7}
\end{equation}%
In the local coordinates we get
\begin{equation}
\alpha (t)=(x(t),y(t)),  \label{4.8}
\end{equation}%
and
\begin{equation}
\overset{.}{\alpha }(t)=\frac{dx}{dt}\frac{\partial }{\partial x}+\frac{dy}{%
dt}\frac{\partial }{\partial y}.  \label{4.9}
\end{equation}%
Now, by means of (\ref{4.7}), from (\ref{4.5}) and (\ref{4.9}), we deduce
the equations so-called \textit{para}-\textit{Hamiltonian equations}
\begin{equation}
\frac{dx}{dt}=\frac{\cosh 2y}{\cosh ^{2}2y-1}\frac{\partial H}{\partial y},%
\frac{dy}{dt}=-\frac{\cosh 2y}{\cosh ^{2}2y-1}\frac{\partial H}{\partial x}.
\label{4.10}
\end{equation}%
In the end, we may say to be \textit{para-mechanical system }$(\mathbf{R}%
^{2},\Phi ,Z_{H})$ triple on para-K\"{a}hlerian space form $(\mathbf{R}%
^{2},g,J).$

\section{Discussion}

From above,\textbf{\ }we understand that Lagrangian and Hamiltonian
formalisms in generalized Classical Mechanics and field theory can be
intrinsically characterized on $(\mathbf{R}^{2},g,J)$ being a model of para-K%
\"{a}hlerian space forms$.$ So, the paths of semispray $\xi $ on $\mathbf{R}%
^{2}$ are the solutions of the Euler-Lagrange equations given by (\ref{3.12}%
) on the mechanical system $(\mathbf{R}^{2},\Phi _{L},\xi )$. Also, the
solutions of the Hamiltonian equations determined by (\ref{4.10}) on the
mechanical system $(\mathbf{R}^{2},\Phi ,Z_{H})$ are the paths of vector
field $Z_{H}$ on $\mathbf{R}^{2}$.

\end{document}